\magnification 1200\baselineskip=20pt
\settabs 10 \columns

\def\B#1{{\bf #1}}
\def\C#1{\centerline{#1}}
\def\CB#1{\C{\B {#1}}}
\def\Sect#1{\noindent\centerline{\bf{#1}}} 



\font\small=cmr7

\parindent=15pt


\CB{Agnes Mary Clerke and Mathematics in the Eleventh {\sl Encyclopedia Britannica\/}}
\C{James S.~Wolper}
\C{Department of Mathematics and Statistics}
\C{Idaho State University}
\C{Pocatello, ID~~83209--8085}
\C{\tt wolpjame@isu.edu}

\bigskip

\noindent
{\small Abstract: Agnes Mary Clerke (1842--1907) was an important
contributor to the discussions of Mathematics in the Eleventh Edition of the 
Encyclopedia Britannica.
This is a survey of the Encyclopedia's treatment of Mathematics 
and an appreciation of her writings.}

\bigskip



\bigskip\Sect{Introduction}

It is difficult to imagine that there was once a time when one
person, or one book, or even one collection of books, ``knew
everything."  E.~T.~Bell [B] wrote that Henri
Poincar\'e did: ``[he] was the last man
to take practically all mathematics, both pure and applied, as his
province."  
The Eleventh Edition of the {\sl Encyclopedia Britannica\/} had the 
audacity to ``presume to define and delineate the entire (and thus complete) opus of 
all human knowledge." [N]  Its 29 volumes appeared
simultaneously in 1910 [EB].  There were gestures toward egalitarianism:
the Encyclopedia also
published various study guides that suggested articles to be read by people
in various walks of life, a kind-of early 20th Century MOOC.
A ``Handy Volume" edition, was sold in the United States
through Sears, Roebuck and Co.

Could its 29 volumes match its contemporary, Poincar\'e?

Encyclopedias and dictionaries are problematic for historical inquiries.
At first glance these works  appear to be secondary sources;
this is especially true of current products.
But encyclopedias are written and edited, and older editions offer
a view of the intellectual and social climates of the time of
production.  These works express a point of view and as such may be
treated as a compendium of primary sources, which will become clear.

Some authors are identified, {\it viz.\/}~Samuel Johnson's
{\sl  A Dictionary of the English Language\/} (1755), but many are anonymous,
or, at best, listed in front matter only.  Many authors were prominent intellectuals,
so offer real insight into the era.  For example, one of 
the {\sl Oxford English Dictionary's\/} anonymous authors was 
J.~R.~R.~Tolkien [W], and the 
article on ``Diagrams" in the Ninth Edition of {\sl Britannica\/}
was written by James Clerk Maxwell.
This level of authority applied to encyclopedic treatments of
Mathematics, as detailed below.  The Eleventh {\sl Britannica\/} and previous 
editions were written by a mix of (almost) anonymous editors and named contributors.

To further illustrate {\sl Britannica\/}'s insight into intellectual life, consider that
it did not shy away from editorializing then (the current online edition
is far more bland), and many academic disagreements were continued
on its pages.  Those familiar with the Wade--Giles system of transliterating
Chinese are probably unaware that there was an intense although largely
one--sided rivalry between Rev.~Jame Legge, the first Oxford 
Professor of Chinese, and Herbert Allen Giles, the first Cambridge
Professor of Chinese.  Giles, the aggressor in this dispute,
wrote the description of Chinese characters for {\sl Britannica\/};
the article contains several pointed references to a competing theory
of the the construction of the characters, dismissing it as lacking
in evidence, although Legge was not named.  

More relevant to {\sl Britannica\/}'s discussions of
Mathematics is the mention in the unsigned article on the Bernoulli family that
they were ``driven from [Holland] by the 
oppressive government of Spain."  {\sl Britannica\/}
was a weapon of the Empire, and Spain was a recent enemy (the Battle of 
Trafalgar was in 1805).  The ``entire opus of human knowledge"
valued British knowledge above all else, perhaps
grudgingly admitting that other European
countries had made contributions.

Who wrote about Mathematics for the Encyclopedia?
It is almost axiomatic that most authors in the Eleventh {\sl Britannica\/} and
earlier editions were men:
women in England did not even get the right to vote until 1928.  But
Agnes Mary Clerke contributed at least 29 signed articles or sections of articles
to the work (one suspects that there were other contributions that went
unattributed; it will be seen later that
Clerke was comfortable with anonymity).  The list of contributors to 
Volume 10 of the Handy Volume Issue, in which Clerke's
articles on Lagrange and Laplace appear, has four female
names; there are 210 listed contributors.

Clerke is widely celebrated for her contributions to Astronomy (see [Br1], 
[Br2], [M]); [Br2], written by an astronomer, details these contributions
but includes very little mathematics.
Clerke's mathematical contributions and the place of 
Mathematics in the Eleventh {\sl Britannica\/} are less well-known.

\bigskip\Sect{Clerke and the Encyclopedia}

The biography [Br2] of Clerke by M.~T.~B\"urke shows that Clerke's
exposure to mathematics and science, an important component
of her mathematical contributions, began as a child.  Clerke's
father had studied Mathematics, and her 
brother Aubrey also studied Mathematics at Trinity College, Dublin.  
Her childhood took place in a household
oriented toward mathematics and science.
She also 
knew Latin, Greek, German, French, Italian, and Spanish, in
addition to, of course, English, and reviewed books in those
languages for the {\sl Edinburgh Review\/}.  The Review was edited by Henry Reeve,
who was a contributor to the Ninth {\sl Britannica\/}; Br\"uck ([Br1], [Br2]) makes
the reasonable speculation
that it was through Reeve that Clerke was invited to contribute.

Articles in the {\sl Edinburgh Review\/} were unsigned, so did not
contribute to Clerke's public renown.   Her first was ostensibly a 
review of Leopoldo Franchetti e Sidney Sonnino, {\it
La Sicilia nel 1876\/} but given the running
title ``Brigandage in Sicily."  She also cited
{\it Relazione della Giunta per l'Inchiesta sulle condizioni della Sicilia
nominata secondo il disposto dell' Articolo 2 della Legge 3 Luglio 1875\/} [Report
of the commission to investigate the condition of Sicily created following the 
promulgation of Article 2 of the law of July 3, 1875].  Clerke's review is
quite opinionated and includes quotes from Cicero and Shakespeare.

Most notable about Clerke's first publication was how modern it seems.
Her discussion of Sicilian society necessarily includes material about
the various {\it mafia\/} and their local control, and it reads like 
Mario Puzo's {\sl The Godfather\/} [P] and movies like {\sl Prizzi's Honor\/} [Pr].
This speaks to her ability to find the essence of a subject, that is, what is likely
to be regarded as important at a later date.

\bigskip\Sect{Mathematics in the Encyclopedia}

Searching the Eleventh to understand its treatment of 
Mathematics, 
one sees several different styles of article.  Some are straightforward,
but some are surprising.  

Olaus Magnus Friedrich Henrici's article on ``Calculating Machines" 
reads like a more modern textbook, including the now--common usage of the 
royal ``we," as in ``We now give an example$\ldots$."  As to depth,
it includes derivations of the integrals involved.

The article ``Groups, Theory Of" provides another illustration of how deep  the 
Encyclopedia's content could be.  It is wedged between ``Ground Squirrel" and ``Grouse,"
which some might find amusing, but it was written by William Burnside, who
is well-known for his contributions to the theory of groups in general.  
The choice of topics is enigmatic.
There is no concept in the article corresponding to an abstract group; groups
are described as transformations of sets of objects.
Surprisingly, it begins with a development of the theory of
continuous groups, which Burnside attributed to Sophus Lie,  including derivations and
a detailed discussion of the one-parameter subgroups generated by an 
infinitesimal operation, although he did not use
the word ``exponential."  Then
Burnside outlined the concept of factor groups, not using the word
``normal," and immediately proceeded to discuss composition series.  One
might judge it too terse to learn from.  It is interesting to note
that, unlike Giles, Burnside was generous in his acknowledgement
of the contributions of others.

Another interesting example is the article ``Fourier's {\it [sic\/]} Series" by Ernest William
Hobson.  One supposes that Burnside, who made major contributions
to the Representation Theory of Groups, could have written this
article.  Hobson, who was on the faculty at Cambridge,
began by introducing the concept of pointwise convergence
of a series of functions, then went on to introduce uniform convergence.
The corresponding article in the contemporary online 
{\sl Britannica\/} does use the word ``convergence" but
does not draw the distinction between pointwise and uniform convergence. 

Evidently Britannica authors approached
writing with different ideas about the interest and experience of the
readers; this is not unusual in Mathematics even today.  Agnes Mary Clerke
had a clearer concept of her readers in mind.

\bigskip\Sect{Clerke's Influence and the Modern Encyclopedia\/}

Clerke wrote the articles on Lagrange and Laplace for the
Encyclopedia.  These were fairly long and included hints about
personalities (Bell [B] did this, too).
Clerke's articles on Lagrange and Laplace focus on their contributions to
Astronomy and Celestial Mechanics.  She did not include Lagrange's 
Theorem on the order of
subgroups.  (Burnside derived the theorem in the article on groups, but
did not use Lagrange's name.)  While she discussed Lagrange's
contributions to the Calculus of Variations, she
did not discuss the technique of Lagrange multipliers.

Clerke's article on Laplace was mixed in its judgments.  Her judgment of Laplace's popular
{\it Exposition du syst\`eme du monde\/} (Paris, 1796) was ``The style is lucid and
masterly, and the summary of astronomical history with which it terminates has been reckoned
[no attribution; one suspects this is the judgment of
Clerke herself] one of the masterpieces of the language."

But she wrote of Laplace's entry into politics ``genius degraded to 
servility for the sake of a ribbon and a title."  She also quotes Napoleon's
judgment ``He brought into the administration the spirit of the
infinitesimals."

One must bear in mind that the Encyclopedia, like the OED
and, conjecturally, the Bible (see, {\it e.g.,\/} the article
on Julius Wellhausen in the Eleventh), is a large work compiled from the
contributions of many authors, and as such is subject to wide variability
in style, depth, and outlook.  The editors presumably insured that the
outlook was consistently pro--British, but as seen with
Burnside's appropriate citation of Lie this could not have been 
excessively rigid.  The editors' oversight could not ensure consistency
of depth or style.  (The variation in style and the personal nature of 
many articles is one of the sources of pleasure in reading {\sl Britannica\/}.)

Most of Clerke's articles were written for the Ninth Edition of
{\sl Britannica\/} and were reprinted in the Eleventh Edition;
she died before the Eleventh appeared.
She was the subject of an unsigned article in the Eleventh {\sl Britannica\/}. 
The article seems confused.
``English astronomer and scientific writer," it begins, before going on to
say Clerke was ``not a practical astronomer in the ordinary sense."  
An editor concerned with length could have omitted mention
of ``astronomer" rather asserting her to be such then deprecating
that assertion.  This theme is common in contemporaneous
assessments of Clerke's work, many of which taint the
objective judgment of her work by noting that 
it is that of a ``lady." [Br2]

Some modern writers devalue her work, at least in the
realm of Mathematics.  Br\"uck [(Br1], [Br2]) speculates about Clerke's
article on Laplace, asserting that it 
borrowed heavily from her brother's prize--winning essay on 
theoretical mechanics because ``the material is treated mathematically
rather than descriptively."  This judgment seems questionable on
several
grounds.  First, the version of Clerke's article in the Eleventh Edition
(unlike many others in the Encyclopedia) contains no mathematical
formulas, and does not even mention partial differential equations or integral
transforms, so the assertion that
the treatment is mathematical rather than descriptive does not hold
up to scrutiny.   The article does go into some depth in discussing 
Laplace's development of probability and of generating functions, 
at about the same of detail as the article on Lagrange goes into his
work in the theory of equations.  With no other articles about
pure mathematics by her with which to compare one should probably withhold
judgment in this matter.

Second, make note that
an encyclopedia article is presumed to offer a summary of the
its subject (although this is not always the case), and as such it would be expected that 
several experts might make suggestions or collaborate.
The article on Astronomy, jointly signed by Clerke and
Simon Newcomb, is in two sections; Clerke wrote the historical
section, Newcomb the scientific.  In either case her familiarity with
primary sources is clear.

Another resolution of these difficulties lies in a different appreciation of 
Clerke as primarily a writer rather than primarily a scientist who made
no claim to be a mathematician.  But she would have been familiar with much
Mathematics from the practice of astronomy, as made clear in her discussion of
John Herschel's career at Cambridge in [C1].  For example, this
work tells the story of Herschel's efforts (along with Peacock and
Babbage) to introduce so--called continental algebraic
techniques into the practice of Mathematics at Cambridge University, especially
in regards to notation.  In fact her story about Herschel's characterization
of the persistence of Newton's dot notation for derivatives as
``dot$\cdot$age" is mentioned in Struik's {\sl Concise History of 
Mathematics\/} [S] although without attribution.

Contrast Clerke's treatment of mathematical subjects to that
of E.~T.~Bell in {\sl Men of Mathematics\/} [B].  Bell included 
many formulas and sketches of proofs, which were appropriate in his
work but many might feel inappropriate in a work for the general public.
Bell had the advantage of being a single author with a clear
concept of his audience.
[By the way, Bell was confused about the Herschels and their r\^ole; his
index lists two mentions of the elder Herschel, but in fact one of them
applies to the younger.]

Nobody accused Clerke of being ``not a member of a crime syndicate in the
ordinary sense" after her article on Sicilian criminal activity.
Clerke was a prolific writer, with
an emphasis on astronomy.  But in addition to Sicily and her
many contributions to the {\sl Edinburgh Review\/} she wrote about Homer [C4].
Even in that work her emphasis on astronomy comes to the fore.
She included a chapter on ``Homer's Astronomy," and the chapter
``The Dog in Homer," while beginning with an amusing if not perhaps
over--written appreciation of the animal's good and bad points, suddenly
switches to a discussion of Sirius, the Dog Star, before returning to
close textual analysis.  The question she addressed is whether Homer loved
or hated dogs given the difference in their treatment in the {\sl Iliad\/}
and the {\sl Odyssey\/}, and this led to speculation about the 
actual authors of these poems as well as other Greek
verse of that era.

This discussion illustrates how Clerke could be obscure.  She
wrote ``$\ldots$the bard of Odysseus has long ceased to possess
[a name].  His only appellation must remain for all time that of
his hero in the Cyclops' cave," but she did not go on to remind the reader
that Odysseus had told Cyclops that his name was ``Nobody." [H]

Clerke wrote essays about many topics, and
her {\sl Britannica\/} biographer praised her ``remarkable skill
in collating, interpreting and summarizing."   Few prolific essayists
write about Mathematics.  This, as much as being a successful
female intellectual in an anti--feminist age, distinguishes her.


\bigskip\Sect{References}

\medskip [Note: articles in the Eleventh Edition of {\sl
Britannica\/} are cited as such in the text and not listed here.  These
were all found in  {\sl The Encyclop{\ae}dia Britannica\/}, a Dictionary
of Arts, Sciences, Literature, and General Information.  Eleventh Edition.
Handy Volume Issue.  New York: The Encyclop{\ae}dia Britannica Company (1910).]

\medskip
\frenchspacing

\item{[B]} E.~T.~Bell, {\sl Men of Mathematics\/}. NY: Simon and Schuster (n.d.)

\item{[Br1]} Br\"uck, M.~T., ``Agnes Marie Clerke, a Chronicler of
Astronomy," {\sl Q.~J.~R.~astr.~Soc.\/} (1994),
{\bf 35}, 59--79.

\item{[Br2]} Br\"uck, M.~T., {\sl Agnes Mary Clerke and the Rise of Astrophysics\/},
Cambridge Univ.~Press (2002).

\item{[C1]} Clerke, Agnes Mary, {\sl The Herschels and Modern Astronomy\/},
New York: MacMillan \& Co. (1895)

\item{[C2]} Clerke, Agnes Mary, {\sl A Popular History of Astronomy during the Nineteenth 
Century,\/} London: A.~and C.~Black (1902).

\item{[C3]} Clerke, Agnes Mary, ``Brigandage in Sicily," {\sl Edinburgh 
Review\/},  (1877), {\bf CXLV}, 251--60.

\item{[C4]} Clerke, Agnes Mary, {\sl Familiar studies in Homer\/}.
 London, New York, Longmans, Green (1892).

\item {[EB]} ``Eleventh Edition And Its Supplements," https://www.britannica.com/topic/

Encyclopaedia-Britannica-English-language-reference-work/Eleventh-edition-and-

its-supplements,

accessed 12 February 2019.

\item{[H]} Homer ({\it attr.\/}), {\sl The Odyssey\/}, Book 9.

\item{[M]} ``Agnes Marie Clerke" at MacTutor History 
of Mathematics archive, 

http://www-groups.dcs.st-and.ac.uk/history/Biographies/Clerke.html,

accessed 15 February 2019.

\item {[N]} Unsigned, ``Eleventh Edition," {\sl The New Yorker\/},
March 2, 1981.

\item {[P]} Puzo, Mario, {\sl The Godfather\/}.  NY: G.~P.~Putnam's Sons (1969).

\item {[Pr]} {\sl Prizzi's Honor\/}, John Huston, dir. (1985).

\item{[S]} Struik, Dirk, {\sl A Concise History of Mathematics\/}, NY: Dover.

\item{[W]} Winchester, Simon, {\sl The Meaning of Everything\/}.

\bye